\documentstyle[12pt]{article}
\begin{document}
\centerline{\bf QUOJECTIONS WITHOUT BANACH SUBSPACES}
\centerline{\bf M.I.Ostrovskii\footnote[1]{Supported by T\"UB\.ITAK}}
\centerline{Mathematics Department, Bo\v gazi\c{c}i University}
\centerline{80815 Bebek, Istanbul, TURKEY}
\centerline{and}
\centerline{Mathematical Division, Institute for Low Temperature Physics}
\centerline{47 Lenin avenue, 310164 Kharkov, UKRAINE}

e-mail: mostrovskii@ilt.kharkov.ua


{\bf Abstract} A quojection (projective limit of Banach spaces
with surjective linking mappings) without infinite dimensional
Banach subspaces is constructed. This results answers a question
posed by G.Metafune and V.B.Moscatelli.\footnote[2]{1991
Mathematics Subject Classification 46A13, 46B20}
\bigskip

Recall that a Fr\'echet space
$F$ is called a {\it quojection} if there exists a sequence
$\{F_n\}_{n=1}^\infty$ of Banach spaces and a sequence of linear
continuous surjective maps $R_n:F_{n+1}\to F_n\ (n\in{\bf N})$ such
that $F$ is isomorphic to the projective limit of the system
$\{R_n\}_{n=1}^\infty$. A
quojection is called {\it twisted} if it is not isomorphic to
a countable product of Banach spaces. The first example of
twisted quojection was constructed by V.B.Moscatelli [M]. This
example and its different generalizations are important tools
in Fr\'echet space theory and its applications (see [MM2]).
G.Metafune and V.B.Moscatelli [MM1, p. 149] (see also [MM2, p. 248])
asked whether there exist twisted quojections without infinite
dimensional Banach subspaces? Such quojections may be considered
as ``very twisted''. The purpose of the present paper is to
construct such quojections. This result may be also considered
as a strengthening of a result due to J.C.Diaz [D].
Our source for Banach space
theory is [LT].

{\bf Theorem.} {\it There exists a twisted quojection without infinite
dimensional Banach subspaces.}

Proof. It is known [MM1, p. 149] that in order to prove Theorem 1
it is sufficient to find a sequence $\{F_n\}_{n=1}^\infty$ of
Banach spaces satisfying the following conditions.

1. For every $n\in{\bf N}$ there exists a surjective linear
continuous mapping $R_n:F_{n+1}\to F_n.$

2. For every $j,k\in{\bf N},\ j\ne k$ the spaces $F_j$ and $F_k$
do not have isomorphic infinite-dimensional subspaces.

In order to do this we use the following version of the construction due to
S.F.Bellenot [B].

Let $1<q<\infty$, $Z$ be a separable Banach space and let
$\{Z_n\}_{n=1}^\infty$ be an increasing sequence of finite-dimensional
subspaces of $Z$ such that cl$(\cup_{n=1}^\infty Z_n)=Z$. Let
$Z_0=\{0\}$.

Let $\{z_i\}_{i=0}^\infty$ be a finitely non-zero sequence with
$z_i\in Z_i$. We define the norm $||\cdot||_{J,q}$ of this sequence
by
$$2||\{z_i\}_{i=0}^\infty||_{J,q}^q=\sup
\{\sum_{i=1}^{k-1}||z_{\eta(i)}-z_{\eta(i+1)}||^q+
||z_{\eta(k)}||^q\}^{1/q},$$
where the sup is over all sequences $\{\eta(i)\}_{i=1}^k$ of integers
with $0\le\eta(1)<\dots<\eta(k).$ The completion of the space of
all finitely non-zero sequences $\{z_i\}_{i=0}^\infty$ under the
norm $||\cdot||_{J,q}$ we shall denote by $J(q,\{Z_n\}_{n=0}^\infty)$.
We shall denote by $J(q,Z)$ the class of the spaces obtained on this
way for fixed separable Banach space $Z$ and all possible sequences
$\{Z_n\}_{n=0}^\infty$.

{\bf Lemma 1.} {\it The quotient $X^{**}/X$ is isometric to $Z$
for every $X\in J(q,Z)$.}

In [B, p. 99] this statement was proved for $q=2$ but it is easy
to verify that the given proof works for every
$1<q<\infty$.

It is easy to verify that the following generalizations of
Lemmas 3 and 5 from [O] are valid.

{\bf Lemma 2.} {\it Let $Z$ be a separable Banach space,
$X\in J(q,Z)$. Then every closed infinite dimensional
subspace of $X$ contains subspaces isomorphic to} $l_q$.

{\bf Lemma 3.} {\it Let $Z$ be a separable Banach space,
$X\in J(q,Z)$ and $r>q$. Then $X^{**}$ does not contain
subspaces isomorphic to} $l_r$.

Using these lemmas we are able to prove

{\bf Lemma 4.} {\it Let a separable Banach space $Z$ be such
that every of its closed infinite dimensional subspaces contains
a subspace isomorphic to $l_r$, $r>q>1$
and $X\in J(q,Z)$. Then every closed infinite dimensional subspace
of $X^{**}$ contains a subspace isomorphic to} $l_q$.

Proof. By Lemma 1 $X^{**}/X$ is isometric to $Z$. Therefore
by the well-known facts (see [LT], Chapter 2) about the
structure of spaces $l_p$ and Lemma 3 it follows that the
quotient mapping $X^{**}\to X^{**}/X$ is strictly singular.
It follows that every infinite dimensional subspace of
$X^{**}$ contains an infinite dimensional subspace isomophic
to a subspace of $X$. Applying  Lemma 2 we obtain the desired
result. $\Box$

Let a sequence $\{q(n)\}_{n=1}^\infty$ of real numbers be such that
$1<q(n)<\infty$ and $q(n+1)<q(n)$. Let us introduce Banach spaces
$\{F_n\}_{n=1}^\infty$ in the following way: $F_1=l_{q(1)}$ and
for $n\ge 2$ let $F_n=X^{**}$, where $X$ is arbitrarily chosen
space from $J(q(n), F_{n-1})$. By Lemma 1 there exist quotient
mappings $R_n:F_{n+1}\to F_n\ (n\in{\bf N})$. Using Lemma 4 and
induction on $n$ we prove that every closed infinite dimensional
subspace of $F_n$ contains a subspace isomorphic to $l_{q(n)}$.
By known facts about spaces $l_p$ (see [LT], Chapter 2) it follows
that Condition 2 is also satisfied. The theorem is proved. $\Box$
\newpage
\centerline{\bf References}

[B] S.F.Bellenot,
The $J$-sum of Banach spaces,
J. Funct. Anal.
48
(1982), pp. 95--106

[D] J.C.Diaz, An example of Fr\'echet space, not Montel, without
infinite dimensional normable subspaces, Proc. Amer. Math. Soc.
96 (1986), p. 721

[LT] J. Lindenstrauss and L. Tzafriri,
Classical Banach
spaces I, Sequence spaces, Springer-Verlag, Berlin, 1977

[MM1] G.Metafune and V.B.Moscatelli,
On twisted Fr\'echet and (LB)-spaces,
Proc. Amer. Math. Soc.
108
(1990), pp. 145--150

[MM2] G.Metafune and V.B.Moscatelli,
Quojections and prequojections,
Advances in the Theory of Fr\'echet spaces
(edited by T.Terzio\v glu), Kluwer Academic Publishers, Dordrecht,
1989, pp. 235--254

[M] V.B.Moscatelli, Fr\'echet spaces without continuous norms
and without bases, Bull. London Math. Soc. 12 (1980), pp. 63--66

[O] M.I.Ostrovskii, Subspaces containing biorthogonal functionals of
bases of different types (Russian), Teor. Funktsii,  Funktsional.  Anal. i
   Prilozhen. 57 (1992), pp. 115--127
\end{document}